\documentclass[11pt,reqno]{amsart}
\usepackage{enumerate, latexsym, amsmath, amsfonts, amssymb, amsthm, color}
\def\pmod #1{\ ({\rm{mod}}\ #1)}
\def\Z{\Bbb Z}

\def\l{\left}
\def\r{\right}
\def\bg{\bigg}
\def\({\bg(}
\def\){\bg)}
\def\t{\text}
\def\f{\frac}
\def\mo{{\rm{mod}\ }}

\def\ls{\leqslant}
\def\gs{\geqslant}

\def\sm{\setminus}
\def\bi{\binom}

\def\ve{\varepsilon}

\def\eq{\equiv}

\def\FF#1#2#3{{}_2F_1\bigg(\bmatrix{#1}\\{#2}\end{bmatrix}\bigg|#3\bigg)}
\def\Proof{\noindent{\it Proof}}
\def\Ack{\medskip\noindent {\bf Acknowledgments}}
\theoremstyle{plain}
\newtheorem{theorem}{Theorem}

\newtheorem{lemma}{Lemma}
\newtheorem{corollary}{Corollary}
\newtheorem{conjecture}{Conjecture}
\theoremstyle{definition}

\theoremstyle{remark}
\newtheorem{remark}{Remark}
\renewcommand{\theequation}{\arabic{section}.\arabic{equation}}
\renewcommand{\thetheorem}{\arabic{section}.\arabic{theorem}}
\renewcommand{\thelemma}{\arabic{section}.\arabic{lemma}}
\renewcommand{\thecorollary}{\arabic{section}.\arabic{corollary}}
\renewcommand{\theconjecture}{\arabic{conjecture}}
\renewcommand{\theremark}{\arabic{remark}}

\vspace{4mm}
\begin{document}
\renewcommand{\baselinestretch}{1.3}
\renewcommand{\arraystretch}{1.3}

\hbox{Nanjing Univ. J. Math. Biquarterly 36 (2019), no.\,2, 134-155.}
\medskip

\title
[{Some new problems in additive combinatorics}]
{Some new problems in additive combinatorics}

\author
[Zhi-Wei Sun] {Zhi-Wei Sun}

\address {Department of Mathematics, Nanjing
University, Nanjing 210093, People's Republic of China}
\email{zwsun@nju.edu.cn}

\thanks{
\newline \indent Supported by the National Natural Science Foundation (grant 11571162)
 of China.}

 \begin{abstract} In this paper we investigate some new problems in additive combinatorics.
 Our problems mainly involve permutations (or circular permutations) $n$ distinct numbers
 (or elements of an additive abelian group)  $a_1,\ldots,a_n$ with adjacent sums $a_i+a_{i+1}$ (or differences $a_i-a_{i+1}$)
 pairwise distinct. For an odd prime power $q=2n+1>13$ with $q\not=25$, we show that
 there is a circular permutation $(a_1,\ldots,a_n)$ of the elements of $S=\{a^2:\ a\in\mathbb F_q\setminus\{0\}\}$
 such that $\{a_1+a_2,\ldots,a_{n-1}+a_n,a_n+a_1\}=S$, where $\mathbb F_q$ denotes the field of order $q$.
 For any finite subset $A$ of an additive torsion-free abelian group $G$ with $|A|=n>3$, we prove that
 there is a numbering $a_1,\ldots,a_n$ of the elements of $A$ such that
 $$a_1+2a_2,\ a_2+2a_3,\ \ldots,\ a_{n-1}+2a_n,\ a_n+2a_1$$
 are pairwise distinct. We also pose 30 open conjectures for further research.
\end{abstract}

\maketitle

\renewcommand{\baselinestretch}{1.3}
\renewcommand{\arraystretch}{1.3}

\vskip 0.5cm
 \noindent{\Large\bf 1\quad Introduction}
 \vskip 0.5cm

\renewcommand{\theequation}{1.\arabic{equation}}
\renewcommand{\thetheorem}{1.\arabic{theorem}}
\renewcommand{\thecorollary}{1.\arabic{corollary}}
\renewcommand{\thelemma}{1.\arabic{lemma}}
\renewcommand{\theremark}{1.\arabic{remark}}
\setcounter{lemma}{0} \setcounter{theorem}{0}
\setcounter{corollary}{0} \setcounter{remark}{0}
\setcounter{equation}{0} \setcounter{conjecture}{0}

Additive combinatorics is an active field involving both number theory and combinatorics.
 For an excellent introduction to problems and results in this fascinating field, one may consult Tao and Vu \cite{TV}.
 See also Alon \cite{A} for a useful tool called Combinatorial Nullstellensatz.
 In this paper we study some new problems in additive combinatorics, they involve some special kinds of permutations or circular permutations.

 Now we present our basic results.

\begin{theorem}\label{Th1.1} Let $a_1,\ldots,a_n$ be a monotonic sequence of $n>1$ distinct real
 numbers. Then there is a permutation $(b_1,\ldots,b_n)$ of
 $a_1,\ldots,a_n$ with $b_1=a_1$ such that
 $$|b_1-b_2|,\ |b_2-b_3|,\ \ldots,\ |b_{n-1}-b_n|$$
 are pairwise distinct.
 \end{theorem}
 \begin{remark}\label{Rm1.1}\rm  Theorem 1.1 is the starting point of our topics in this paper.
 \end{remark}

\begin{corollary}\label{Cor1.1} There is a circular permutation
$(q_1,\ldots,q_n)$ of the first $n>2$ primes $p_1,\ldots,p_n$ with
$q_1=p_1=2$ and $q_n=p_n$ such that the $n$ distances
$$|q_1-q_2|,\ |q_2-q_3|,\ \ldots,\ |q_{n-1}-q_n|,\ |q_n-q_1|$$
are pairwise distinct.
\end{corollary}
\Proof. By
Theorem 1.1, there is a permutation $(-q_n,-q_{n-1},\ldots,-q_2)$ of
$-p_n,-p_{n-1},\ldots,-p_2$ with $q_n=p_n$ such that
$|-q_n+q_{n-1}|,\ldots,|-q_3+q_2|$ are pairwise distinct. Set
$q_1=p_1=2$. Then the circular permutation $(q_1,q_2,\ldots,q_n)$ of $p_1,p_2,\ldots,p_n$
and meets our requirement since $|q_1-q_2|=q_2-2$ and $|q_n-q_1|=p_n-2$ are both odd while
those $q_i-q_{i+1}\ (1<i<n)$ are even. \qed

\begin{theorem}\label{Th1.2} {\rm (i)} For any integer $n>3$, there is a circular permutation $(i_0,\ldots,i_n)$ of $0,\ldots,n$
with $i_0=0$ and $i_n=n$ such that all the $n+1$ adjacent differences
$i_0-i_1,i_1-i_2,\ldots,i_{n-1}-i_n,i_n-i_0$ are pairwise distinct.

{\rm (ii)} An integer $n>1$ is even if and only if
there is a permutation $(i_1,\ldots,i_n)$ of $1,\ldots,n$ with
$$i_1-i_2,\ i_2-i_3,\ \ldots,\ i_{n-1}-i_n$$
pairwise distinct modulo $n$.
\end{theorem}
\begin{remark}\label{Rem1.2} \rm In contrast with Theorem 1.2(i), for any $n>2$ distinct integers $a_1<\ldots<a_n$ we clearly have
$$a_1+a_2<a_2+a_3<\ldots<a_{n-1}+a_n.$$
On Sept. 13, 2013 the author asked his students the following question:
When $a_n+a_1=a_i+a_{i+1}$ for some $1\ls i<n$, how to construct a suitable permutation $b_1,\ldots,b_n$ of $a_1,\ldots,a_n$
such that $b_1+b_2,b_2+b_3,\ldots,b_{n-1}+b_n,b_n+b_1$ are pairwise distinct? The author's PhD student Dianwang Hu suggested that
it suffices to take $(b_1,\ldots,b_n)=(a_1,\ldots,a_i,a_{i+2},a_{i+1},a_{i+3},\ldots,a_n)$.
But this does not work for $i=n-2$. If $i>2$, then the permutation
$(b_1,\ldots,b_n)=(a_1,\ldots,a_{i-2},a_i,a_{i-1},a_{i+1},a_{i+2},\ldots,a_n)$
meets the requirement. The case $n=3$ is trivial. For $n=4$, the permutation
$(a_1,a_2,a_4,a_3)$ works for our purpose since $a_1+a_2<a_3+a_1<a_2+a_4<a_4+a_3$.
\end{remark}

\begin{theorem}\label{Th1.3} For any $n>3$ distinct elements $a_1,a_2,\ldots,a_n$ of a torsion-free abelian group $G$,
there is a circular permutation $(b_1,\ldots,b_n)$ of $a_1,\ldots,a_n$ such that all the $n$ sums
$$b_1+2b_2,\ b_2+2b_3,\ \ldots,\ b_{n-1}+2b_n,\ b_n+2b_1$$
are pairwise distinct.
\end{theorem}
\begin{remark}\rm For any circular permutation $(a_1,a_2,a_3)$ of $0,1,2$, the three numbers
$a_1+2a_2,a_2+2a_3,a_3+2a_1$ cannot be pairwise distinct.
\end{remark}

\begin{theorem}\label{Th1.4} For any odd prime power $n>2$, there are integers $a_1,a_2,\ldots,a_{\varphi(n)}$ such that
both $\{a_1,\ldots,a_{\varphi(n)}\}$ and
$$\{a_1-a_2,\ a_2-a_3,\ \ldots,\ a_{\varphi(n)-1}-a_{\varphi(n)},\ a_{\varphi(n)}-a_1\}$$
are reduced systems of residues modulo $n$, where $\varphi$ is Euler's totient function.
\end{theorem}
\begin{remark} \rm  We conjecture that Theorem \ref{Th1.4} remains valid if we just let $n>2$ be odd.
\end{remark}

\begin{theorem}\label{Th1.5} Let $\mathbb F_q$ be a finite field with $q=2n+1>5$ an odd prime power,
and set $S=\{a^2:\ a\in \mathbb F_q\sm\{0\}\}$.

{\rm (i)} If $q\not\in\{9,25\}$, then there is a circular permutation
$(a_1,\ldots,a_n)$ of all the
$n$ elements of $S$ such that
$$\{a_1-a_2,\ a_2-a_3,\ \ldots,\ a_{n-1}-a_n,\ a_n-a_1\}=S.$$

{\rm (ii)} Suppose that $q>13$ and $q\not=25$. Then
there is a circular permutation
$(b_1,\ldots,b_n)$ of all the
$n$ elements of $S$ such that
$$\{b_1+b_2,\ b_2+b_3,\ \ldots,\ b_{n-1}+b_n,\ b_n+b_1\}=S.$$
\end{theorem}
\begin{remark} \rm  In the initial version of this paper posted to arXiv in 2013, the author posed the following conjecture
weaker than Theorem 1.5 which was later confirmed by N. Alon and J. Bourgain \cite{AB}: For any prime $p=2n+1>5$,
there is a circular permutation $(a_1,\ldots,a_n)$ of the $(p-1)/2=n$
quadratic residues modulo $p$ such that all the $n$ adjacent differences $a_1-a_2,a_2-a_3,\ldots,a_{n-1}-a_n,a_n-a_1$ are quadratic residues
 modulo $p$. Also, for any prime $p=2n+1>13$, there is a circular permutation $(b_1,\ldots,b_n)$ of the $(p-1)/2=n$
quadratic residues modulo $p$ such that all the $n$ adjacent sums $b_1+b_2,b_2+b_3,\ldots,b_{n-1}+b_n,b_n+b_1$ are quadratic residues
 modulo $p$.
\end{remark}
\medskip

We are going to prove Theorems 1.1-1.5 in the next section, and pose sixteen conjectures concerning permutations in Section 3.

The author \cite{S17}  proved that for any integer $m>4$ there is a number $n\in\Z^+=\{1,2,3,\ldots\}$ such that $\pi(mn)=m+n,$
where $\pi(x)$ denotes the number of primes not exceeding $x$. Sun  (cf.
\cite[Conjecture 4.4]{S17} and \cite[A247824]{S}) also conjectured that for any $m\in\Z^+$ there is a positive integer $n$ such that $m+n$ divides $p_m+p_n$, where $p_k$ denotes the $k$-th prime.
This has been verified for all $m=1,\ldots,10^5$. (The reader may consult \cite{CP} for basic knowledge about asymptotic behaviors of $\pi(x)$ and $p_n$.)
With this background, we pose in Section 4 fourteen conjectures involving addition and divisibility.

We have posted to OEIS some sequences (cf. \cite{S}) related to our conjectures in Sections 3 and 4.

\medskip

\vskip 0.5cm
 \noindent{\Large\bf 2\quad Proofs of Theorems 1.1-1.5}
 \vskip 0.5cm

\renewcommand{\theequation}{2.\arabic{equation}}
\renewcommand{\thetheorem}{2.\arabic{theorem}}
\renewcommand{\thecorollary}{2.\arabic{corollary}}
\renewcommand{\thelemma}{2.\arabic{lemma}}
\renewcommand{\theremark}{2.\arabic{remark}}

 \setcounter{lemma}{0} \setcounter{theorem}{0}
\setcounter{corollary}{0}
\setcounter{remark}{0}
\setcounter{equation}{0}
\setcounter{conjecture}{0}

\medskip
\noindent{\it Proof of Theorem 1.1}. If $a_1>a_2>\ldots>a_n$, then $-a_1<-a_2<\ldots<-a_n$. So we may
assume that $a_1<a_2<\ldots<a_n$ without loss of generality.

 If $n=2k$ is even, then the permutation
 $$(b_1,\ldots,b_n)=(a_1,a_{2k},a_2,a_{2k-1},\ldots,a_{k-1},a_{k+2},a_k,a_{k+1})$$
 meets our purpose since
 $$a_{2k}-a_1>a_{2k}-a_2>a_{2k-1}-a_2>\ldots>a_{k+2}-a_{k-1}>a_{k+2}-a_k>a_{k+1}-a_k.$$
 When $n=2k-1$ is odd, the permutation
$$(b_1,\ldots,b_n)=(a_1,a_{2k-1},a_2,a_{2k-2},\ldots,a_{k-1},a_{k+1},a_k)$$
 meets the requirement since
 $$a_{2k-1}-a_1>a_{2k-1}-a_2>a_{2k-2}-a_2>\ldots>a_{k+1}-a_{k-1}>a_{k+1}-a_k.$$
This concludes the proof. \qed

\medskip
\noindent{\it Proof of Theorem 1.2}. (i) We first assume that $n=2k$ is even. If $k$ is even, then the circular permutation
$$(i_0,\ldots,i_n)=(0,2k-1,1,2k-2,2,\ldots,k+1,k-1,k,2k)$$
meets the requirement since
$$-(2k-1),\ 2k-2,\ -(2k-3),\ 2k-4,\ \ldots,\ 2,\ -1,\ -k,\ 2k$$
are pairwise distinct. If $k$ is odd, then it suffices to choose the circular permutation
$$(i_0,\ldots,i_n)=(0,1,2k-1,2,2k-2,\ldots,k-1,k+1,k,2k)$$
since
$$-1,\ -(2k-2),\ 2k-3,\ -(2k-4),\ \ldots,\ -2,\ 1,\ -k,\ 2k$$
are pairwise distinct.

 Now we handle the case $n=2k+1\eq1\ (\mo\ 2)$. If $k$ is even, then the circular permutation
$$(i_0,\ldots,i_n)=(0,2k,1,2k-1,2,2k-2,\ldots,k-1,k+1,k,2k+1)$$
meets the requirement since
$$-2k,\ 2k-1,\ -(2k-2),\ 2k-3,\ -(2k-4),\ \ldots,\ -2,\ 1,\ -(k+1),\ 2k+1$$
are pairwise distinct.
If $k$ is odd, then it suffices to choose the circular permutation
$$(i_0,\ldots,i_n)=(0,k,k+2,k+1,k-1,k+3,k-2,k+4,k-3,\ldots,2k-1,2,2k,1,2k+1)$$
since
$$-k,\ -2,\ 1,\ 2,\ -4,\ 5,\ \ldots,\ -(2k-2),\ 2k-1,\ -2k,\ 2k+1$$
are pairwise distinct.

(ii) Suppose that $i_1,\ldots,i_n$ is a permutation of $1,\ldots,n$ with
the $n-1$ integers $i_k-i_{k+1}\ (0<k<n)$ pairwise distinct modulo $n$. Then
$$\{i_k-i_{k+1}\ \mo\ n: \ k=1,\ldots,n-1\}=\{r\ \mo\ n:\ r=1,\ldots,n-1\}$$
and also
$$\{i_{k+1}-i_k\ \mo\ n: \ k=1,\ldots,n-1\}=\{r\ \mo\ n:\ r=1,\ldots,n-1\}.$$
Therefore
$$\sum_{k=1}^{n-1}(i_k-i_{k+1})\eq\sum_{r=1}^{n-1}r\eq \sum_{k=1}^{n-1}(i_{k+1}-i_k)\pmod n$$
and hence $n\mid 2(i_1-i_n)$ which implies that $n$ is even.

Now assume $n=2m$ with $m\in\Z^+$. Then
$$(i_1,\ldots,i_n)=(m,m-1,m+1,m-2,m+2,\ldots,2,2m-2,1,2m-1,2m)$$
is a permutation of $1,\ldots,n$ with the required property.

In view of the above, we have completed the proof of Theorem 1.2. \qed

\medskip
\noindent{\it Proof of Theorem 1.3}. The subgroup of $G$ generated by $a_1,\ldots,a_n$ is a finitely generated torsion-free abelian group.
So we may simply assume that $G=\Z^r$ for some positive integer $r$ without any loss of generality.
It is well known that there is a linear ordering $\ls$ on $G=\Z^r$ such that for any $a,b,c\in G$
if $a<b$ then $-b<-a$ and $a+c<b+c$. For convenience we suppose that $a_1<a_2<\ldots<a_n$ without any loss of generality.

Clearly $a_1+2a_2<a_2+2a_3<\ldots<a_{n-1}+2a_n$.
Thus the permutation $(b_1,\ldots,b_n)=(a_1,\ldots,a_n)$ meets the requirement if
$a_n+2a_1\not=a_i+2a_{i+1}$ for all $i=1,\ldots,n-1$.

Below we assume that $a_n+2a_1=a_i+2a_{i+1}$ for some $0<i<n$. Note that $1\ls i\ls n-2$ since
$a_{n-1}+2a_n-(a_n+2a_1)=a_{n-1}+a_n-2a_1>0$.
\medskip

{\it Case} 1. $i=1$.

In this case, $a_n+2a_1=a_1+2a_2$ and hence $a_1+a_3<a_1+a_n=2a_2$. The permutation
$(b_1,\ldots,b_n)=(a_1,a_3,a_2,a_4,\ldots,a_n)$ meets our purpose since
$$a_n+2a_1=a_1+2a_2<a_1+2a_3<a_3+2a_2<a_2+2a_4<\ldots<a_{n-1}+2a_n.$$

{\it Case} 2. $i>1$ and $n=4$.

In this case, $a_4+2a_1=a_2+2a_3$ and we may take the permutation
$(b_1,b_2,b_3,b_4)=(a_2,a_1,a_3,a_4)$ since
$$a_2+2a_1<a_1+2a_3<a_2+2a_3=a_4+2a_1<a_4+2a_2<a_3+2a_4.$$

{\it Case} 3. $i\gs2$, $n\gs 5$, and  $a_{i-1},a_i,a_{i+1}$ don't form an AP (arithmetic progression).

In this case, the permutation
$$(b_1,\ldots,b_n)=(a_1,\ldots,a_{i-1},a_{i+1},a_i,a_{i+2},\ldots,a_n)$$
works for our purpose since
\begin{align*}&\min\{a_{i-1}+2a_{i+1},a_{i+1}+2a_i\}
\\<&\max\{a_{i-1}+2a_{i+1},a_{i+1}+2a_i\}<a_i+2a_{i+1}=a_n+2a_1
\\<&a_i+2a_{i+2}<\ldots<a_{n-1}+2a_n.
\end{align*}

{\it Case} 4.  $2\ls i<n-2$ and  $a_i-a_{i-1}=a_{i+1}-a_i\not=a_{i+2}-a_{i+1}$.

In this case, the permutation
$$(b_1,\ldots,b_n)=(a_1,\ldots,a_{i-1},a_i,a_{i+2},a_{i+1},a_{i+3},\ldots,a_n)$$
works for our purpose since
\begin{align*} a_{i-1}+2a_i<&a_i+2a_{i+1}=a_n+2a_1
\\<&\min\{a_i+2a_{i+2},a_{i+2}+2a_{i+1}\}<\max\{a_i+2a_{i+2},a_{i+2}+2a_{i+1}\}
\\<&a_{i+1}+2a_{i+3}<\ldots<a_{n-1}+2a_n.
\end{align*}

{\it Case} 5. $2\ls i<n-2$, and $a_{i-1},a_i,a_{i+1},a_{i+2}$ form an AP.

In this case, the permutation
$$(b_1,\ldots,b_n)=(a_1,\ldots,a_{i-1},a_{i+2},a_{i+1},a_i,a_{i+3},\ldots,a_n)$$
works for our purpose since
\begin{align*} a_{i+1}+2a_i<&a_i+2a_{i+1}=a_n+2a_1
\\<&a_{i-1}+2a_{i+2}\ (\t{since}\ a_i-a_{i-1}=a_{i+2}-a_{i+1}<2(a_{i+2}-a_{i+1}))
\\<&a_{i+2}+2a_{i+1}=a_i+2a_{i+2}<a_i+2a_{i+3}<\ldots<a_{n-1}+2a_n.
\end{align*}

{\it Case} 6. $i=n-2\gs3$ and $a_{i+1}-a_i=a_i-a_{i-1}\not=a_{i-1}-a_{i-2}$.

In this case, the permutation
$$(b_1,\ldots,b_n)=(a_1,\ldots,a_{i-2},a_i,a_{i-1},a_{i+1},a_{i+2})$$
works for our purpose since
\begin{align*} \min\{a_{i-2}+2a_i,a_i+2a_{i-1}\}
<&\max\{a_{i-2}+2a_i,a_i+2a_{i-1}\}<a_{i-1}+2a_i
\\<&a_{i-1}+2a_{i+1}<a_i+2a_{i+1}=a_n+2a_1
\\<&a_{i+1}+2a_{i+2}=a_{n-1}+2a_n.
\end{align*}

{\it Case} 7. $i=n-2\gs3$, and $a_{i-2},a_{i-1},a_i,a_{i+1}$ form an AP.

In this case, the permutation
$$(b_1,\ldots,b_n)=(a_1,\ldots,a_{i-2},a_{i+1},a_i,a_{i-1},a_{i+2})$$
works for our purpose since
\begin{align*} a_{i-2}+2a_{i-1}<&a_i+2a_{i-1}
\\<&a_{i-2}+2a_{i+1}\ (\t{since}\ a_i-a_{i-2}=a_{i+1}-a_{i-1}<2(a_{i+1}-a_{i-1}))
\\<&a_{i+1}+2a_i=a_{i-1}+2a_{i+1}<a_i+2a_{i+1}=a_n+2a_1
\\<&a_{i-1}+2a_{i+2}=a_{n-3}+2a_n.
\end{align*}

Combining the above we have finished the proof of Theorem 1.3.

\medskip
\noindent{\it Proof of Theorem 1.4}. Write $n=p^a$
with $p$ an odd prime and $a$ a positive integer. Take a primitive root $g$ modulo $n$.
Note that $g\not\eq1\pmod{p}$ and $g^{\varphi(n)+1}\eq g\pmod n$.
Clearly, both $\{g^i:\ i=1,\ldots,\varphi(n)\}$ and
$$\{g^i-g^{i+1}=g^i(1-g):\ i=1,\ldots,\varphi(n)\}$$
are reduced systems of residues modulo $n$.  So it suffices to take $a_i=g^i$ for $i=1,\ldots,\varphi(n)$. \qed

\begin{lemma}\label{Lem-g} Let $q$ be an odd prime power and set
$S=\{a^2:\ a\in \mathbb F_q\sm\{0\}\}.$

{\rm (i)} The field $\mathbb F_q$ has a primitive element $g$ with $g^2-1\in S$ if and only if $q\not\in\{3,5,9,25\}$.

{\rm (ii)} The field $\mathbb F_q$ has a primitive element $g$ with $g^2+1\in S$ if and only if $q\not\in\{3,5,7,9,13,25\}$.
\end{lemma}
\Proof. For an odd prime $p$ let $G_p$ be the set of those integers $g\in\{\pm1,\ldots,\pm(p-1)/2\}$
which are primitive roots modulo $p$. Then
$$G_3=\{-1\},\ G_5=\{\pm1\},\ G_7=\{-2,3\},\ G_{13}=\{\pm2,\pm6\}$$
and
$$G_{61}=\{\pm2,\pm6,\pm7,\pm10,\pm17,\pm18,\pm26,\pm30\}.$$
Note that
$$\l(\f{3^2-1}7\r)=1,\ \l(\f{2^2-1}{13}\r)=1,\ \t{and}\ \l(\f{2^2+1}{61}\r)=1=\l(\f{7^2-1}{61}\r),$$
where $(-)$ denotes the Legendre symbol.
Thus, it is easy to see that both parts of Lemma \ref{Lem-g} hold for $q\in\{3,5,7,13,61\}$.

The field $\mathbb F_9$ can be viewed as
$$\mathbb F_3[x]/(x^2+1)=\{ax+b\ \mo\ x^2+1:\ a,b\in \mathbb F_3=\{0,\pm1\}\},$$
and it has four primitive elements: $\pm x\pm1\ \mo\ x^2+1$. In the ring $\mathbb F_{3}[x]$ we have the congruence
$$(x\pm1)^2=x^2\mp x+1\eq\mp x\pmod{x^2+1}.$$
So, for any primitive element $g$ of $\mathbb F_9$ neither $g^2+1$ nor $g^2-1$ belongs to $\{a^2:\ a\in \mathbb F_9\}$.

Similarly, $\mathbb F_{25}$ can be viewed as
$$\mathbb F_5[x]/(x^2-2)=\{ax+b\ \mo\ x^2-2:\ a,b\in \mathbb F_5=\{0,\pm1,\pm2\}\},$$
and it has eight primitive elements:
$$\pm x\pm2\ \mo\ x^2-2,\ \ \t{and}\ \ \pm2x\pm1 \ \mo\ x^2-2.$$
In the ring $\mathbb F_{5}[x]$ we have the congruences
$$(x\pm2)^2=x^2\mp x-1\eq1\mp x\pmod{x^2-2}\ \t{and}\ (2x\pm1)^2=-x^2\mp x+1\eq-1\mp x\pmod{x^2-2}.$$
So, for any primitive element $g$ of $\mathbb F_{25}$ neither $g^2+1$ nor $g^2-1$ belongs to $\{a^2:\ a\in \mathbb F_{25}\}$.

The field $\mathbb F_{121}$ can be viewed as
$$\mathbb F_{11}[x]/(x^2-2)=\{ax+b\ \mo\ x^2-2:\ a,b\in \mathbb F_{11}=\{0,\pm1,\ldots,\pm5\}\},$$
and its primitive elements are the polynomials
$$\pm x\pm2,\ \pm2x\pm2,\ \pm2x\pm4,\ \pm3x\pm3,\ \pm3x\pm5,\ \pm4x\pm4,\ \pm5x\pm1,\ \pm5x\pm5$$
modulo $x^2-2$.
In the ring $\mathbb F_{11}[x]$ we have the congruences
$$(2x+4)^2+1\eq (4x+2)^2\pmod{x^2-2}\ \t{and}\ (x+2)^2-1\eq(3x-3)^2\pmod{x^2-2}.$$
So, there are primitive elements $g_1$ and $g_2$ of $\mathbb F_{121}$
with $g_1^2+1,g_2^2-1\in\{a^2:\ a\in \mathbb F_{121}\sm\{0\}\}$.

Now we consider the remaining case: $q>13$ and $q\not\in\{25,61,121\}$.
Let $\ve\in\{\pm1\}$. By \cite[Corollary 2]{COS}, there exists a primitive element $g$ of $\mathbb F_q$
with $g+\ve g^{-1}$ also primitive. Let $n=(q-1)/2$. As
$g^n=-1=(g+\ve g^{-1})^n$, we have $(g^2+\ve)^n=g^n(g+\ve g^{-1})^n=1$ and hence $g^2+\ve\in S$.

In view of the above, we have completed the proof of Lemma \ref{Lem-g}. \qed

\medskip
\noindent{\it Proof of Theorem 1.5}. (i) As $q>5$ and $q\not\in\{9,25\}$, by Lemma \ref{Lem-g}(i),
there is a primitive element $g_0$ of $\mathbb F_q$ with $g_0^2-1\in S$.
Note that $g=g_0^{-1}$ is also a primitive element of $\mathbb F_q$ and $1-g^2\in S$.
Clearly, $\{g^{2i}:\ i=1,\ldots,n\}=S$ and $g^{2(n+1)}=g^2$.
Observe that the $n$ elements
$g^{2i}-\ve g^{2(i+1)}=g^{2i}(1-g^2)\ (i=1,\ldots,n)$
are pairwise distinct and they all belong to $S$.
This proves part (i) of Theorem 1.5.

(ii) As $q>13$ and $q\not=25$, by Lemma \ref{Lem-g}(ii),
there is a primitive element $g$ of $\mathbb F_q$ with $g^2+1\in S$.
Note that $g^2=g^{2(n+1)}$ and
$$\{g^{2i}:\ i=1,\ldots,n\}=S=\{g^{2i}(1+g^2)=g^{2i}+g^{2(i+1)}:\ i=1,\ldots,n\}.$$
So, part (ii) of Theorem 1.5 also holds. \qed
\medskip

\vskip 0.5cm
 \noindent{\Large\bf 3\quad Some Conjectures concerning Permutations}
 \vskip 0.5cm

\renewcommand{\theequation}{3.\arabic{equation}}
\renewcommand{\thetheorem}{3.\arabic{theorem}}
\renewcommand{\thecorollary}{3.\arabic{corollary}}
\renewcommand{\thelemma}{3.\arabic{lemma}}
\renewcommand{\theconjecture}{3.\arabic{conjecture}}
\renewcommand{\theremark}{3.\arabic{remark}}

 \setcounter{lemma}{0} \setcounter{theorem}{0}
\setcounter{corollary}{0}
\setcounter{remark}{0}
\setcounter{equation}{0}
\setcounter{conjecture}{0}

\begin{conjecture}\label{Conj3.1} {\rm (2013-09-01)} Let $a_1,a_2,\ldots,a_n$ be $n>1$ distinct real numbers. Then
there is a permutation $(b_1,\ldots,b_n)$ of $a_1,\ldots,a_n$ with $b_1=a_1$ such that the $n-1$ numbers
$$|b_1-b_2|,\ |b_2-b_3|,\ \ldots,\ |b_{n-1}-b_n|$$
are pairwise distinct.
\end{conjecture}
\begin{remark} \rm By Theorem 1.1, Conjecture 3.1 holds when $a_1$ is the least element or the largest element of $\{a_1,\ldots,a_n\}$. After learning this conjecture from the initial version
of this paper posted to arXiv in 2013, F. Monopoli \cite{M} managed to prove that Conjecture 3.1 holds if the set $A = \{a_1,a_2,\ldots,a_n\}$ forms an arithmetic progression.
\end{remark}

\begin{conjecture}\label{Conj3.3} {\rm (2013-09-02)} Let $a_1,\ldots,a_n$ be $n$ distinct elements of a finite additive abelian group $G$.
Suppose that $n\nmid|G|$, or $n$ is even and the Sylow $2$-subgroup of $G$ is cyclic.
Then there exists a permutation $(b_1,\ldots,b_n)$ of $a_1,\ldots,a_n$
with $b_1=a_1$ such that the $n-1$ elements $b_i-b_{i+1}\ (0<i<n)$ are pairwise distinct.
\end{conjecture}
\begin{remark}\label{Rem3.3} \rm By Theorem 1.2(ii), Conjecture \ref{Conj3.3} holds when $\{a_1,\ldots,a_n\}=G=\Z/n\Z$ with $n$ even.
For the Klein quaternion group
$$G=\Z/2\Z\oplus\Z/2\Z=\{(0,0),(0,1),(1,0),(1,1)\},$$
if $\{a_1,a_2,a_3,a_4\}=G$ then we have $a_1-a_2=a_3-a_4$.
\end{remark}

A subset $A$ of a set $S$ with $|A|=n\in\Z^+$ is called an $n$-subset of $S$.

\begin{conjecture}\label{Conj3.4} {\rm (2013-09-03)} Let $n>2$ be an integer, and let $A$ be an $n$-subset of a finite additive abelian group $G$ with $|G|$ odd.

{\rm (i)} There always exists a numbering $a_1,a_2,\ldots,a_n$ of all the $n$ elements of $A$ such that the $n$ sums
$$a_1+a_2,\ a_2+a_3,\ \ldots,\ a_{n-1}+a_n,\ a_n+a_1$$
are pairwise distinct.

{\rm (ii)} If $n\nmid |G|$, then there is a numbering $a_1,a_2,\ldots,a_n$ of all the $n$ elements of $A$ such that the $n$ differences
$$a_1-a_2,\ a_2-a_3,\ \ldots,\ a_{n-1}-a_n,\ a_n-a_1$$
are pairwise distinct.
\end{conjecture}
\begin{remark}\label{Rem3.4} \rm  A conjecture of H. S. Snevily \cite{Sn} states that for any two $n$-subsets $A$ and $B$ of an additive abelian group of odd order, we may write $A=\{a_1,\ldots,a_n\}$ and $B=\{b_1,\ldots,b_n\}$ so that the $n$
sums $a_1+b_1,\ldots,a_n+b_n$ are pairwise distinct. This was proved by Arsovski \cite{A} in 2009. Note that Conjecture \ref{Conj3.4}(i)
is stronger than Snevily's conjecture in the case $A=B$. In the spirit of Remark 1.2, the assertion in Conjecture \ref{Conj3.4}(i) holds if $A$ is an $n$-subset of a torsion-free abelian group with $n>2$.
See \cite[A228762]{S} for some data related to Conjecture \ref{Conj3.4}(ii).
\end{remark}

\begin{conjecture}\label{Conj3.5}
Let $A$ be a finite subset of an additive abelian group $G$ with $|A|=n>3$.

{\rm (i)} {\rm (2013-09-20)} If $G$ is finite with $|G|\not\eq0\ (\mo\ 3)$, then there is a numbering $a_1,\ldots,a_n$
of all the elements of $A$ such that the $n$ sums
$$a_1+2a_2,\ a_2+2a_3,\ \ldots,\ a_{n-1}+2a_n,\ a_n+2a_1$$
are pairwise distinct.

{\rm (ii)} {\rm (2017-12-11)} Let $B$ be any $n$-subset of $G$. Then we can write $A=\{a_1,\ldots,a_n\}$
and $B=\{b_1,\ldots,b_n\}$ so that either the $n$ sums
$$a_1+2b_1,\ a_2+2b_2,\ \ldots,\ a_{n-1}+2b_{n-1},\ a_n+2b_n$$
are pairwise distinct, or the $n$ sums
$$2a_1+b_1,\ 2a_2+b_2,\ \ldots,\ 2a_{n-1}+b_{n-1},\ 2a_n+b_n$$
are pairwise distinct.
\end{conjecture}
\begin{remark}\label{Rem3.5}\rm  (i) When $A=\{a_1,\ldots,a_n\}$ forms an abelian group
$G$ of the form $(\Z/3\Z)^r$, the $n$ elements
$$a_1+2a_2=a_1-a_2,\ a_2+2a_3=a_2-a_3,\ \ldots,\ a_{n-1}+2a_n=a_{n-1}-a_n,\ a_n+2a_1=a_n-a_1$$
cannot be pairwise distinct (otherwise $0\not\in\{a_1-a_2,a_2-a_3,\ldots,a_n-a_1\}=G$).

(ii) The author has proved part (ii) in the case $A=B$ for $n=4$.
When $G$ is an abelian group of odd order, Conjecture \ref{Conj3.5}(ii) is equivalent to Snevily's conjecture confirmed by Arsovski \cite{A} since $a=2(\f{|G|+1}2a)$ for all $a\in G$.
\end{remark}
\medskip

\begin{conjecture}\label{Conj3.7} {\rm (joint with Qing-Hu Hou)} {\rm  (i) (2013-09-05)}
For any finite field $\mathbb F_q$ with $q>7$, there is a numbering
$a_1,\ldots,a_q$ of all the elements of $\Bbb F_q$ such that all the $q$ sums
$$a_1+a_2,\ a_2+a_3,\ \ldots,\ a_{q-1}+a_q,\ a_q+a_1$$
are primitive elements of $\Bbb F_q$.

{\rm (ii) (2013-09-07)} Let $p=2n+1$ be an odd prime. If $p>19$, then there is a circular permutation $(i_1,\ldots,i_n)$ of $1,\ldots,n$ such that
all the $n$ adjacent sums $i_1+i_2,i_2+i_3,\ldots,i_{n-1}+i_n,i_n+i_1$ are primitive roots modulo $p$. When $p>13$,
there is a circular permutation $(i_1,\ldots,i_n)$ of $1,\ldots,n$ such that
all the $n$ adjacent differences $i_1-i_2,i_2-i_3,\ldots,i_{n-1}-i_n,i_n-i_1$ are primitive roots modulo $p$.
\end{conjecture}
\begin{remark}\label{Conj3.7}\rm (a) We have verified part (i) for primes $q$ below $545$, and part (ii) for primes $p$ below $545$. For the field $$\mathbb F_9=\mathbb F_3[x]/(x^2+1)=\{\overline{ax+b}:\ a,b\in\mathbb F_3=\{0,\pm1\}\}$$
with $\overline{ax+b}$ the residue class of $ax+b$ modulo $x^2+1$,
the four primitive elements are $\overline{\pm x\pm1}$. For the circular permutation
$$(\bar 0,\overline{-x+1},\overline{1},\overline{-x},\overline{-x-1},\overline{-1},\overline{x-1},
\overline{x},\overline{x+1})$$
of all the elements of $\mathbb F_9$, all the $9$ adjacent sums are primitive elements of $\mathbb F_9$.
For the circular permutation
$$(a_1,a_2,\ldots,a_{11})=(0,6,7,1,5,3,10,8,9,4,2)$$
of $0,1,\ldots,10$, the $11$ sums $a_1+a_2,a_2+a_3,\ldots,a_{10}+a_{11},a_{11}+a_1$ are all primitive roots modulo the prime $11$.

(b) If $g$ is a primitive element of the field $\Bbb F_q$ with $q>3$ and $a_i=g^{i-1}$ for all $i=1,\ldots,q-1$, then it is easy to see that
$a_1-a_2,a_2-a_3,\ldots,a_{q-2}-a_{q-1},a_{q-1}-a_1$ are pairwise distinct and that
$a_1+a_2,a_2+a_3,\ldots,a_{q-2}+a_{q-1},a_{q-1}+a_1$ are also pairwise distinct.
\end{remark}

\begin{conjecture}\label{Conj3.8} {\rm (2013-09-11)} Let $p=2n+1$ be an odd prime. If $p>19$, then there is a circular permutation
$(a_1,\ldots,a_n)$ of all the
$(p-1)/2=n$ quadratic residues modulo $p$ such that all the $n$ adjacent sums $a_1+a_2,a_2+a_3,\ldots,a_{n-1}+a_n,a_n+a_1$
are primitive roots modulo $p$. If $p>13$, then there is a circular permutation
$(b_1,\ldots,b_n)$ of all the
$(p-1)/2=n$ quadratic residues modulo $p$ such that all the $n$ differences $b_1-b_2,b_2-b_3,\ldots,b_{n-1}-b_n,b_n-b_1$
are primitive roots modulo $p$.
\end{conjecture}
\begin{remark}\label{Rem3.8}\rm  Compare this conjecture with Theorem 1.5
and Conjecture \ref{Conj3.7}(ii).
\end{remark}

\begin{conjecture}\label{Conj3.9} {\rm (2013-09-15)} Let $p=2n+1>11$ be a prime.

{\rm (i)} There is a circular permutation $(i_1,\ldots,i_n)$ of $1,\ldots,n$
such that all the $n$ numbers $i_1^2+i_2,\ i_2^2+i_3,\ \ldots,\ i_{n-1}^2+i_n,\ i_n^2+i_1$
are quadratic residues modulo $p$. Also, there is a circular permutation $(j_1,\ldots,j_n)$ of  $1,\ldots,n$
such that all the $n$ numbers $j_1^2-j_2,\ j_2^2-j_3,\ \ldots,\ j_{n-1}^2-j_n,\ j_n^2-j_1$
are quadratic residues modulo $p$.

{\rm (ii)} If $p>13$, then there is a circular permutation $(i_1,\ldots,i_n)$ of $1,\ldots,n$
such that all the $n$ numbers
$i_1^2+i_2,\ i_2^2+i_3,\ \ldots,\ i_{n-1}^2+i_n,\ i_n^2+i_1$
are primitive roots modulo $p$. Also, there is a circular permutation $(j_1,\ldots,j_n)$ of  $1,\ldots,n$
such that all the $n$ numbers
$j_1^2-j_2,\ j_2^2-j_3,\ \ldots,\ j_{n-1}^2-j_n,\ j_n^2-j_1$
are primitive roots modulo $p$.
\end{conjecture}
\begin{remark}\label{Rem3.9}\rm For example, $(i_1,\ldots,i_{11})=(1,6,7,11,4,5,3,8,10,9,2)$ is a circular permutation of $1,\ldots,11$ for which
all the sums $i_1^2+i_2,i_2^2+i_3,\ldots,i_{10}^2+i_{11},i_{11}^2+i_1$ are primitive roots modulo $23$. Also,
$(j_1,\ldots,j_{11})=(1,9,7,5,11,10,3,2,6,8,4)$ is a circular permutation of $1,\ldots,11$ for which
all the sums $j_1^2-i_2,j_2^2-i_3,\ldots,j_{10}^2-j_{11},j_{11}^2-i_1$ are primitive roots modulo $23$.
See also \cite[A229141]{S} for related data.
\end{remark}

\begin{conjecture}\label{Conj3.10} {\rm (2013-09-17)} Let $\Bbb F_q$ be a finite field with $q>7$ elements, and let $a_0$ be any element of $\Bbb F_q$. Then
there is a circular permutation $(a_1,\ldots,a_{q-1})$ of all the nonzero elements of $\Bbb F_q$ such that all the $q-1$ elements
$a_0+a_1a_2,a_0+a_2a_3,\ldots,a_0+a_{q-2}a_{q-1},a_0+a_{q-1}a_1$ are primitive elements of the field $\Bbb F_q$.
\end{conjecture}
\begin{remark}\label{Rem3.10}\rm For the circular permutation $(i_1,\ldots,i_{10})=(1,9,2,4,5,8,10,3,6,7)$ of $1,\ldots,10$,
all the $10$ integers $i_1i_2-1,i_2i_3-1,\ldots,i_9i_{10}-1,i_{10}i_1-1$ are primitive roots modulo 11.
\end{remark}

\begin{conjecture}\label{Conj3.11} {\rm (2013-09-07)} For any positive integer $n\not=2,4$, there exists a permutation $i_0,i_1,\ldots,i_n$ of $0,1,\ldots,n$
with $i_0=0$ and $i_n=n$ such that all the $n+1$ adjacent sums
$i_0+i_1,\ i_1+i_2,\ \ldots,\ i_{n-1}+i_n,\ i_n+i_0$
are coprime to both $n-1$ and $n+1$.
\end{conjecture}
\begin{remark}\label{Rem3.11}\rm (i) See \cite[A228886]{S} for related data. Note that there is no circular permutation $i_0,\ldots,i_7$ of $0,\ldots,7$ with $i_0+i_1,i_1+i_2,\ldots,i_6+i_7,i_7+i_0$
all relatively prime to $7\times13-1=90$. We also guess that $n\pm1$ in Conjecture \ref{Conj3.11} can be replaced by $2n\pm1$.

(ii) Now we explain why Conjecture \ref{Conj3.11} holds for any positive odd integer $n$.
If $n\eq1,3\ (\mo\ 6)$, then $n-2$ and $2n-1$ are relatively prime to both $n-1$ and $n+1$, and hence the circular permutation
$$(i_0,\ldots,i_n)=(0,n-2,2,n-4,4,\ldots,1,n-1,n)$$
meets the requirement. If $n\eq3,5\ (\mo\ 6)$, then $n+2$ is relatively prime to both $n-1$ and $n+1$, and hence the circular permutation
$$(i_0,\ldots,i_n)=(0,1,n-1,3,n-3,\ldots,n-2,2,n)$$
suffices for our purpose.
\end{remark}

\begin{conjecture}\label{Conj3.12} {\rm (2013-09-22)} {\rm (i)} Let $A$ be a set of $n>2$ distinct nonzero real numbers. Then
there is a circular permutation $(a_1,a_2,...,a_n)$ of all the elements of $A$ such that the $n$ adjacent sums
$a_1+a_2, a_2+a_3, ..., a_{n-1}+a_n, a_n+a_1$
are pairwise distinct, and that the $n$ adjacent products
$a_1a_2, a_2a_3, ..., a_{n-1}a_n, a_na_1$
are also pairwise distinct, except for the following three cases:

  {\rm (a)} $|A|=4$ and $A$ has the form $\{\pm s,\pm t\}$.

  {\rm (b)} $|A|=5$ and $A$ has the form $\{r,\pm s,\pm t\}$.

  {\rm (c)} $|A|=6$ and $A$ has the form $\{\pm r,\pm s,\pm t\}$.
\medskip

{\rm (ii)} For any set $A$ of $n>3$ distinct nonzero real numbers,
there is a circular permutation $(a_1,a_2,...,a_n)$ of all the elements of $A$ such that the $n$ adjacent differences
   $a_1-a_2, a_2-a_3, ..., a_{n-1}-a_n, a_n-a_1$
are pairwise distinct, and that the $n$ adjacent products
   $a_1a_2, a_2a_3, ..., a_{n-1}a_n, a_na_1$
are also pairwise distinct, unless $|A|=4$ and $A$ has the form $\{\pm s,\pm t\}$.
\end{conjecture}
\begin{remark}\label{Rem3.12} For the set $A=\{1,2,\ldots,n\}$ with $n\in\{3,5,7,\ldots\}$, obviously
$1+2,2+3,\ldots,(n-1)+n,n+1$ are pairwise distinct since $n+1$ is even, and
$1\times 2,2\times 3,\ldots,(n-1)n,n\times 1$ are also pairwise distinct since $n$ is odd.
\end{remark}

\begin{conjecture}\label{Conj3.13} {\rm (2013-09-08)} For any positive integer $n$, there is a circular permutation $(i_0,i_1,\ldots,i_n)$ of $0,1,\ldots,n$
such that all the $n+1$ adjacent sums
$i_0+i_1,\, i_1+i_2,\, \ldots,\, i_{n-1}+i_n,\, i_n+i_0$
belong to the set $\{k\in\Z^+:\ 6k-1\ \t{and}\ 6k+1\ \t{are twin primes}\}$.
\end{conjecture}
\begin{remark}\label{Rem3.13}\rm Clearly this conjecture implies the twin prime conjecture. Qing-Hu Hou has verified this conjecture for all $n\ls 100$.
We also have similar conjectures for cousin primes, sexy primes, and primes of the form $4k-1$ or $4k+1$ or $6k+1$ (cf. \cite[A228917]{S}).
In 1982 A. Filz \cite{F} (see also \cite[p.\,160]{G}) conjectured that for any $n=2,4,6,\ldots$ there is a circular permutation
$i_1,\ldots,i_n$ of $1,\ldots,n$ such that all the $n$ adjacent sums $i_1+i_2,i_2+i_3,\ldots,i_{n-1}+i_n,i_n+i_1$ are prime.
\end{remark}

\begin{conjecture}\label{Conj3.14} {\rm (2013-09-08)} For any integer $n>2$, there exists a circular permutation $(i_0,i_1,\ldots,i_n)$ of $0,1,\ldots,n$
such that all the $n+1$ adjacent sums $i_0+i_1,\, i_1+i_2,\, \ldots,\, i_{n-1}+i_n,\, i_n+i_0$ are of the form $(p+1)/6$,
where $p$ is a Sophie Germain prime.
\end{conjecture}
\begin{remark}\label{Rem3.14}\rm Recall that a prime $p$ with $2p+1$ also prime is called a Sophie Germain prime.
It is conjectured that there are infinitely many Sophie Germain primes.
\end{remark}

\begin{conjecture}\label{Conj3.15} {\rm (i) (2013-09-09)} For any positive integer $n$, there exists a circular permutation $(i_0,i_1,\ldots,i_n)$ of $0,1,\ldots,n$
such that all the $2n+2$ numbers
 $$|i_0\pm i_1|,\ |i_1\pm i_2|,\ \ldots,\ |i_{n-1}\pm i_n|,\ |i_n\pm i_0|$$
 are of the form $(p-1)/2$, where $p$ is an odd prime.

 {\rm (ii) (2013-09-10)} For any positive integer $n\not=2,4$, there exists a circular permutation $(i_0,i_1,\ldots,i_n)$
of $0,1,\ldots,n$ such that all the $n+1$ numbers
 $$|i_0^2-i_1^2|,\ |i_1^2-i_2^2|,\ \ldots,\ |i_{n-1}^2-i_n^2|,\ |i_n^2-i_0^2|$$
 are of the form $(p-1)/2$, where $p$ is an odd prime.
\end{conjecture}
\begin{remark}\label{Rem3.15}\rm  See \cite[A228956 and A229005]{S} for related data. Here are two suitable circular permutations: $(0,1,2,3,5,4,7,8,6,9)$ for $n=9$ in part (i),
and $(i_0,\ldots,i_5)=(0,1,4,5,2,3)$ for $n=5$ in part (ii).
\end{remark}

\begin{conjecture}\label{Conj3.16} {\rm (2013-09-13)} For any positive integer $n\not=4$, there exists a circular permutation $(i_0,i_1,\ldots,i_n)$
of $0,1,\ldots,n$ with $i_0=0$ and $i_n=1$ such that all the $n+1$ numbers
 $i_0^2+i_1,\ i_1^2+i_2,\ \ldots,\ i_{n-1}^2+i_n,\ i_n^2+i_0$
 are of the form $(p-1)/2$, where $p$ is an odd prime.
\end{conjecture}
\begin{remark}\label{Rem3.16}\rm See \cite[A229082]{S} for related data. For $i,j\in\{0,\ldots,n\}$ with $i+j>1$, if $j$ is a multiple of $3$ and $2(i^2+j)+1$ is a prime then
$2i^2+1\not\eq0\ (\mo\ 3)$ and hence $3\mid i$. So, if $(i_0,i_1,\ldots,i_n)$ is a permutation
of $0,1,\ldots,n$ with $i_0=0$ such that all the $n+1$ numbers
$i_0^2+i_1,\, i_1^2+i_2,\, \ldots,\, i_{n-1}^2+i_n,\, i_n^2+i_0$
 are of the form $(p-1)/2$ with $p$ an odd prime, then we must have $i_n=1$ (otherwise, $i_n,i_{n-1},\ldots,i_1$
 are all divisible by 3 which is impossible). To illustrate Conjecture 3.16, we give a desired permutation for $n=20$:
$$(i_0,\ldots,i_{20})=(0,3,12,9,15,18,6,20,19,14,13,4,2,7,16,17,11,10,5,8,1).$$
\end{remark}

\begin{conjecture}\label{Conj3.17} {\rm (2013-09-16)} Let $n$ by any positive integer. Then there exists a circular permutation $(i_0,i_1,\ldots,i_n)$
of $0,1,\ldots,n$ such that all the $n+1$ numbers
 $i_0^2+i_1,\ i_1^2+i_2,\ \ldots,\ i_{n-1}^2+i_n,\ i_n^2+i_0$
 are of the form $(p-1)/4$ with $p$ a prime congruent to $1$ modulo $4$.
 Also, there is a circular permutation $(j_0,j_1,\ldots,j_n)$
of $0,1,\ldots,n$ with $j_0=0$ and $j_n=1$ such that all the $n+1$ numbers
 $j_0^2+j_1,\ j_1^2+j_2,\ \ldots,\ j_{n-1}^2+j_n,\ j_n^2+j_0$
 are of the form $(p+1)/4$ with $p$ a prime congruent to $3$ modulo $4$.
\end{conjecture}
\begin{remark}\label{Rem3.17}\rm See \cite[A227456]{S} for related data. For $i,j\in\{0,\ldots,n\}$ with $i+j>1$, if $j$ is a multiple of $3$ and $4(i^2+j)-1$ is a prime then
$4i^2-1\not\eq0\ (\mo\ 3)$ and hence $3\mid i$. So, if $(j_0,j_1,\ldots,j_n)$ is a permutation
of $0,1,\ldots,n$ with $j_0=0$ such that all the $n+1$ numbers
$j_0^2+j_1,\, j_1^2+j_2,\, \ldots,\, j_{n-1}^2+j_n,\, j_n^2+j_0$
 are of the form $(p+1)/4$ with $p$ a prime congruent to 3 modulo 4, then we must have $j_n=1$ (otherwise, $j_n,j_{n-1},\ldots,j_1$
 are all divisible by 3 which is impossible). To illustrate Conjecture \ref{Conj3.17}, we give two desired permutations for $n=9$:
$$(i_0,\ldots,i_9)=(0,1,2,3,4,6,9,7,8,5)\ \ \t{and}\ \ (j_0,\ldots,j_9)=(0,3,6,9,2,4,5,8,7,1).$$
\end{remark}

\begin{conjecture}\label{Conj3.18} {\rm (2013-09-17)} For any integer $n>5$ with $n\not=13$, there is a circular permutation
$(i_1,i_2,\ldots,i_n)$ of $1,\ldots,n$ such that $i_1i_2-1,i_2i_3-1,\ldots,i_{n-1}i_n-1,i_ni_1-1$ are all prime.
Also, for any positive integer $n>1$ (resp. $n\not=4$), there is a circular permutation
$(i_1,i_2,\ldots,i_n)$ of $1,\ldots,n$ such that $2i_1i_2-1,2i_2i_3-1,\ldots,2i_{n-1}i_n-1,2i_ni_1-1$
(resp. $2i_1i_2+1,2i_2i_3+1,\ldots,2i_{n-1}i_n+1,2i_ni_1+1$) are all prime.
\end{conjecture}
\begin{remark}\label{Rem3.18}\rm See \cite[A229232]{S} for related data. For the circular permutation
$$(i_1,\ldots,i_{23})=(1,6,23,10,9,22,11,18,13,14,21,2,15,4,17,16,5,12,7,20,19,8,3),$$
all the 23 numbers $i_1i_2-1,i_2i_3-1,\ldots,i_{22}i_{23}-1,i_{23}i_1-1$ are primes.
\end{remark}
\medskip

\vskip 0.5cm
 \noindent{\Large\bf 4\quad Conjectures involving Addition and Divisibility}
 \vskip 0.5cm

\renewcommand{\theequation}{4.\arabic{equation}}
\renewcommand{\thetheorem}{4.\arabic{theorem}}
\renewcommand{\thecorollary}{4.\arabic{corollary}}
\renewcommand{\thelemma}{4.\arabic{lemma}}
\renewcommand{\theconjecture}{4.\arabic{conjecture}}
\renewcommand{\theremark}{4.\arabic{remark}}

 \setcounter{lemma}{0} \setcounter{theorem}{0}
\setcounter{corollary}{0}
\setcounter{remark}{0}
\setcounter{equation}{0}
\setcounter{conjecture}{0}

\begin{conjecture}\label{Conj4.1} Let $m$ be any positive integer.

{\rm (i) (2014-09-29)} $m+n$ divides $p_m^2+p_n^2$ for some $n\in\Z^+$.

{\rm (ii) (2014-09-30)} $m+n$ divides $p_{m^2}+p_{n^2}$ for some $n\in\Z^+$. Moreover, for $m>1$ we may require that $n\ls m(m-1)/2$.
\end{conjecture}
\begin{remark}\label{Rem4.1} \rm See \cite[A247975 and A248354]{S} for related data. For $m=4703$, the least $n\in\Z^+$ with $m+n$ dividing $p_m^2+p_n^2$ is $760027770$.
Note also that $2+3$ divides $p_{2^2}+p_{3^2}=7+23=30$.
\end{remark}

\begin{conjecture}\label{Conj4.2} {\rm (2014-09-29)} Let $m\in\Z^+$ and $\ve\in\{\pm1\}$. Then, there is a positive integer $n$ such that
$$p_{mn}\eq\ve\pmod{m+n};$$
moreover, we may require $n\ls m(m-1)/2$ if $m>2$.
\end{conjecture}
\begin{remark}\label{Rem4.2}\rm  For example, $p_{2\times 4}=19\eq1\pmod{2+4}$. See \cite[A248004]{S} for related data.
\end{remark}

\begin{conjecture}\label{Conj4.3} {\rm (2014-09-30)} Let $m\in\Z^+$. Then $m+n$ divides $\pi(m)^2+\pi(n)^2$ for some $n\in\Z^+$.
Also, $m+n$ divides $\pi(m^2)+\pi(n^2)$ for some $n\in\Z^+$.
\end{conjecture}
\begin{remark}\label{Rem4.3}\rm  See \cite[A248044 and A248052]{S} for related data. For example,
$$\pi(5)^2+\pi(12)^2=3^2+5^2=34\eq0\pmod{5+12}$$ and
$$\pi(4^2)+\pi(8^2)=6+18=24\eq0\pmod{4+8}.$$
\end{remark}

\begin{conjecture}\label{Conj4.4} {\rm (2014-09-27)} Let $a$ be any integer with $a\not\eq3\pmod 6$. Define the Lucas sequences $(u_n(a))_{n\gs0}$ and $(v_n(a))_{n\gs0}$
by
$$u_0(a)=0,\ u_1(a)=1,\ u_{n+1}(a)=au_n(a)+u_{n-1}(a)\ (n=1,2,3,\ldots),$$
and
$$v_0(a)=2,\ v_1(a)=a,\ v_{n+1}(a)=av_n(a)+v_{n-1}(a)\ (n=1,2,3,\ldots).$$
Let $m\in\Z^+$. Then there are infinitely many $n\in\Z^+$ with $m+n$ dividing $u_m(a)+u_n(a)$. Also,
there are infinitely many $n\in\Z^+$ with $m+n$ dividing $v_m(a)+v_n(a)$.
\end{conjecture}
\begin{remark}\label{Rem4.4}\rm Note that those $F_n=u_n(1)\ (n=0,1,2,\ldots)$ are the Fibonacci numbers and those $L_n=v_n(1)\ (n=0,1,2,\ldots)$ are the Lucas numbers.
For $m\in\Z^+$, see \cite[A247937 and A247940]{S} for the least $n>m$ with $m+n$ dividing $F_m+F_n$ (or $L_m+L_n$). See also \cite[A248133, A248136, A248137, A248139, A248142]{S} for similar conjectures.
\end{remark}

\begin{conjecture}\label{Conj4.5}  Let $m$ be a positive integer.

{\rm (i)} {\rm (2014-09-29)} $m+n$ divides $\bi{2m}m+\bi{2n}n$ for some $n\in\Z^+$.

{\rm (ii)} {\rm (2014-09-29)} If $m\not=3$, then $m+n$ divides $C_m+C_n$ for some $n\in\Z^+$ with $\gcd(m,n)=1$,
where $C_k$ refers to the Catalan number
$\bi{2k}k/(k+1)=\bi{2k}k-\bi{2k}{k+1}$.

{\rm (iii) (2014-10-01)} For each $m=1,2,3,\ldots$, there is a positive integer $n$ such that $\gcd(m,n)=1$ and $mn\mid C_{m+n}$.
\end{conjecture}
\begin{remark}\label{Rem4.5}\rm See \cite[A248125, A248124 and A248123]{S} for related data. For example, the least $n\in\Z^+$ with $\gcd(9,n)=1$ and $(9+n)\mid (C_9+C_n)$ is $95$; also, $\gcd (4,21)=1$, and $4\times 21$ divides $C_{4\times 21}=4861946401452$.
\end{remark}

\begin{conjecture}\label{Conj4.6} {\rm (2014-09-29)} For any integer $m>6$, there is a positive integer $n$ with $\varphi(m)\varphi(n)\eq0\pmod{m+n}.$
\end{conjecture}
\begin{remark}\label{Rem4.6}\rm  See \cite[A248007]{S} for related data. For example, $\varphi(10)\varphi(14)=4\times 6\eq0\pmod{10+14}$.
\end{remark}

\begin{conjecture}\label{Conj4.7} {\rm (2014-09-29)} Let $m$ be any positive integer. Then
$m+n$ divides $\sigma(mn)$ for some $n\in\Z^+$, where $\sigma(k)$ refers to the sum of all positive divisors of $k$.
\end{conjecture}
\begin{remark}\label{Rem4.7}\rm See \cite[A248008]{S} for related data. For example, $4+6$ divides $\sigma(4\times6)=60$.
\end{remark}

\begin{conjecture}\label{Conj4.8} {\rm (2014-09-29)} Let $m$ be a positive integer.

{\rm (i)} If $m>1$, then $(m+n)\mid\sigma(m)\varphi(n)$ for some $n=1,\ldots,m$.

{\rm (ii)} There is a positive integer $n$ such that $(m+n)\mid\varphi(m)\sigma(n)$. Moreover, we may require $n<2m$ if $m>2$.
\end{conjecture}
\begin{remark}\label{Rem4.8}\rm See \cite[A248029 and A248030]{S} for related data. For example,
$$\sigma(8)\varphi(7)=15\times6=90\eq0\pmod{8+7}$$
and
$$\varphi(2)\sigma(12)=28\eq0\pmod{2+12}.$$
\end{remark}

\begin{conjecture}\label{Conj4.9} {\rm (2014-09-29)} Let $m$ be any positive integer. Then $m+n$ divides $\varphi(m)^2+\varphi(n)^2$ for some $n\in\Z^+$.
Moreover, we may require $n\ls m^2$ except for $m=33$.
\end{conjecture}
\begin{remark}\label{Rem4.9}\rm See \cite[A248035]{S} for related data. For example,
$$\varphi(33)^2+\varphi(1523)^2=20^2+1522^2=2316884\eq0\pmod{33+1523}.$$
\end{remark}

\begin{conjecture}\label{Conj4.10} {\rm (i) (2014-09-29)} For any $m\in\Z^+$, we have $(m+n)\mid(\sigma(m)^2+\sigma(n)^2)$ for some $n\in\Z^+$.

{\rm (ii) (2014-09-30)} For any $m\in\Z^+$, we have $(m+n)\mid(\sigma(m^2)+\sigma(n^2))$ for some $n\in\Z^+$.
\end{conjecture}
\begin{remark}\label{Rem4.10}\rm See \cite[A248036 and A248054]{S} for related data. For example,
$$\sigma(1024)^2+\sigma(2098177)^2=4423875080209\eq0\pmod{1024+2098177}$$
and
$$\sigma(4^2)+\sigma(7^2)=31+57=88\eq0\pmod{4+7}.$$
\end{remark}

\begin{conjecture}\label{Conj4.11} {\rm (2014-10-02)} For $k\in\Z^+$ let $p(k)$ be the number of partitions of $k$ (i.e., unordered ways to write $k$ as a sum of some positive integers with repetitions allowed).

{\rm (i)} For any $m\in\Z^+$, there is a positive integer $n$ such that $m+n$ divides $p(m)+p(n)$.

{\rm (ii)} For any $m\in\Z^+$, there is a positive integer $n$ such that $(m+n)\mid p(mn)$.
\end{conjecture}
\begin{remark}\label{Rem4.11}\rm  See \cite[A248143 and A248144]{S} for related data. For example,
$$p(5)+p(13)=7+101=108\eq0\pmod{5+13}$$
and $$p(6\times 14)=26543660\eq0\pmod{6+14}.$$
\end{remark}

\begin{conjecture}\label{Conj4.12} {\rm (2014-09-30)} Let $m$ be any positive integer. Then $mn$ divides $\varphi(m^2+n^2)$ for some $n\in\Z^+$.
\end{conjecture}
\begin{remark}\label{Rem4.12}\rm  We have verified this for all $m=1,\ldots,1242$, see
\cite[A248058]{S} for related data. For example, for $m=1093$ the least $n\in\Z^+$
with $mn\mid\varphi(m^2+n^2)$ is $57343152$. In fact,
\begin{align*}\varphi(1093^2+57343152^2)=&\varphi(3288237082489753)=3285228630168576
\\=&52416\times1093\times 57343152.
\end{align*}
If $(mq)^2+1$ is prime for some $q\in\Z^+$, then for $n=m^2q$ we have
$$\varphi(m^2+n^2)=\varphi(m^2)\varphi(1+m^2q^2)=\varphi(m^2)m^2q^2\eq0\pmod{mn}.$$
\end{remark}

\begin{conjecture}\label{Conj4.13} {\rm (2014-10-08)} For any $m\in\Z^+$, there is a positive integer $n$ such that $\varphi(m+n)\mid n$.
Moreover, we may require $n\ls m(m-1)$ for $m>1$.
\end{conjecture}
\begin{remark}\label{Rem4.13}\rm See \cite[A248568]{S} for related data. For example, $\varphi(10+40)=20$ divides $40$.
\end{remark}

\begin{conjecture}\label{4.14} {\rm (2014-10-05)} Let $m$ be any positive integer.  Then $p_{m+n}-p_n$ divides $m+n$ for some $n\in\Z^+$.
Also, $p_{m+n}-p_n$ divides $n$ for some $n\in\Z^+$.
\end{conjecture}
\begin{remark}\label{Rem3.4}\rm  See \cite[A248366 and A248369]{S} for related data. For example,
$p_{5+175}-p_{175}=1069-1039=30$ divides $5+175=180$ and $p_{7+80}-p_{80}=449-409=40$ divides $80$.
\end{remark}

\Ack. The initial version of this paper was posted to arXiv in Sept. 2013 with the ID {\tt arXiv:1309.1679}.
The author would like to thank Prof. Noga Alon for helpful comments, and Prof. Qing-Hu Hou for checking many of the author's conjectures
via a computer.


\begin{thebibliography}{99}

\renewcommand{\arraystretch}{1.2}


\bibitem{A} N. Alon, {\it Combinatorial Nullstellensatz}, Combin. Probab. Comput. {\bf 8} (1999), 7--29.

\bibitem{AB}  N. Alon and J. Bourgain, {\it Additive patterns in multiplicative subgroups},
preprint, 2013.

\bibitem{Ar}  B. Arsovski, {\it  A proof of Snevily's conjecture},
Israel J. Math. {\bf 182} (2011), 505--508.

\bibitem{COS} S. D. Cohen, T. Oliveira e Silva and N. Sutherland,
{\it Linear combinations of primitive elements of a finite field},
Finite Fields Appl. {\bf 51} (2018), 388--406.

\bibitem{CP} R. Crandall and C. Pomerance, Prime Numbers: A Computational Perspective, Springer, New York, 2001.

\bibitem{F}  A. Filz, {\it  Problem 1046}, J. Recreational Math.  {\bf 14} (1982), 64; {\bf 15} (1983), 71.

\bibitem{G}  R. K. Guy, Unsolved Problems in Number Theory, 3rd Edition, Springer, New York, 2004.

\bibitem{M} F. Monopoli, {\it Absolute differences along Hamiltonian paths}, Electron. J. Combin. {\bf 22} (2015), no. 3, \#P3.20, 1--8.

\bibitem{Sn}  H. S. Snevily, {\it  The Cayley addition table of $\Z_n$},
Amer. Math. Monthly  {\bf 106} (1999), 584--585.

\bibitem{S}  Z.-W. Sun, Sequences A227456, A228762, A228886, A228917, A228956, A229005, A229038, A229082, A229141, A229232,
     A247824, A247937, A247940, A247975, A248004, A248007, A248008, A248029, A248030, A248035, A248036, A248044, A248052, A248054, A248058,
 A248123, A248124, A248125, A248133, A248136, A248137, A248139, A248142, A248143, A248144, A248354, A248366, A248369, A248568
    in OEIS (On-Line Encyclopedia of Integer Sequences), {\tt http://oeis.org}

\bibitem{S17} Z.-W. Sun, {\it A new theorem on the prime-counting function}, Ramanujan J.
{\bf 42} (2017), 59--67.

\bibitem{TV}  T. Tao and V. H. Vu, Additive Combinatorics,
Cambridge Univ. Press, Cambridge, 2006.

\end{thebibliography}
\end{document}